\documentclass[oneside]{amsart}
\usepackage{amsmath,amsthm}

\numberwithin{equation}{section}

\newtheorem{theorem}{Theorem}
\newtheorem{lemma}{Lemma}[section]
\newtheorem{prop}[lemma]{Proposition}
\newtheorem{cor}{Corollary}

\theoremstyle{definition}
\newtheorem{definition}[lemma]{Definition}
\newtheorem{example}[lemma]{Example}

\theoremstyle{remark}
\newtheorem{remark}[lemma]{Remark}

\DeclareMathOperator{\vol}{vol}
\DeclareMathOperator{\Vol}{Vol}
\DeclareMathOperator{\grad}{grad}

\newcommand{\R}{\mathbb R}
\renewcommand{\phi}{\varphi}
\newcommand{\pd}{\partial}
\newcommand{\ep}{\varepsilon}
\newcommand{\U}{\mathcal U}

\def\be{\begin{equation}}
\def\ee{\end{equation}}

\begin{document}

\title[Local monotonicity of volume]
{Local monotonicity of Riemannian and Finsler volume
with respect to boundary distances}

\author{Sergei Ivanov}
\thanks{Supported by RFBR grants 09-01-12130-ofi-m and 11-01-00302-a}
\address{St.~Petersburg department of V.A.Steklov Institute of Mathematics
of the Russian Academy of Sciences,
191023, Fontanka 27, Saint Petersburg, Russia}
\email{svivanov@pdmi.ras.ru}

\begin{abstract}
We show that the volume of a simple Riemannian metric
on $D^n$ is locally monotone with respect to its
boundary distance function. Namely if $g$ is a simple
metric on $D^n$ and $g'$ is sufficiently close to $g$
and induces boundary distances greater or equal
to those of $g$, then $\vol(D^n,g')\ge \vol(D^n,g)$.
Furthermore, the same holds for Finsler metrics and the
Holmes--Thompson definition of volume.
As an application, we give a new proof of
injectivity of the geodesic ray transform
for a simple Finsler metric.
\end{abstract}

\maketitle

\section{Introduction}

A Riemannian metric $g$ on the $n$-dimensional
disc $D=D^n$ is called \textit{simple} if
the boundary $\pd D$ is strictly convex
with respect to $g$ (that is, its second
fundamental form is positive definite)
and all geodesics in $(D,g)$ are minimizing
and have no conjugate points
(or, equivalently, every pair of points
in $D$ is connected by a unique geodesic which
varies smoothly with the endpoints).
Note that this property
persists under $C^\infty$-small perturbations of the metric.

For a Riemannian metric $g$ on $D$, we denote by
$d_g$ the distance function on $D\times D$
induced by~$g$.
The \textit{boundary distance function} of $g$,
denoted by $bd_g$, is the restriction of $d_g$
to $\pd D\times\pd D$.
That is, $bd_g(x,y)$ is the length of
a $g$-shortest path in $D$ between boundary points
$x$ and~$y$. If $g$ is simple, this shortest path
is the (unique) $g$-geodesic connecting $x$ and $y$.


It is well-known that the volume of a simple metric $g$
is determined by the function $bd_g$
via an explicit formula involving
boundary distances and their derivatives
(cf.\ \cite{Santalo}, \cite{gromov-frm}, \cite{Croke91}).
It is natural to expect that this formula is monotone
with respect to $bd_g$, that is, if another metric
$g'$ satisfies $bd_{g'}\ge bd_g$ pointwise, then
$\vol(D,g')\ge\vol(D,g)$.
However the formula itself is not monotone
if arbitrary functions are allowed in place
of $bd_g$ (see Example \ref{x-2d}).
On the other hand, functions that can be realized
by boundary distances of simple
metrics are rather special, and it might
be the case that the volume formula is
monotone within this class of functions.

This question is a variant of
the minimal filling conjecture (see \cite{BI10}, \cite{I10})
which asserts that a simple Riemannian metric $g$
has the least volume among all metrics whose
boundary distance functions majorize that of~$g$.
This conjecture is related to Michel's boundary
rigidity conjecture~\cite{Michel} about unique
determination of a simple metric by its
boundary distance function.

The minimal filling conjecture is essentially about
finding metrics realizing filling volumes, see \cite{gromov-frm}.
It has been confirmed in a number of special cases.
In dimension~2 the conjecture is proved 
for any simple metric $g$ and any competing metric $g'$ on $D^2$
(\cite{I02}, see also \cite{I09-pre} for the Finslerian case
and \cite{pestov-uhlmann} for boundary rigidity).
However the general case
where the competing metric $g'$ can be on a surface of arbitrary genus
remains open, even for $(D^2,g)$ 
isometric to a subset of the standard hemisphere.
(The latter special case is equivalent to Gromov's Filling Area Conjecture.)
In higher dimensions, filling minimality is established
for flat metrics by Gromov \cite{gromov-frm},
for regions in negatively curved symmetric spaces
by Besson, Courtois and Gallot \cite{BCG},
and for metrics $g$ sufficiently close to a flat or hyperbolic metric
by Burago and Ivanov \cite{BI10,BI-hyperb}.
These results come with equality case analysis
that yields boundary rigidity of the respective metrics.
Croke and Kleiner \cite{CK98} proved boundary rigidity
of some product metrics using a weaker form of volume minimality.
Croke, Dairbekov and Sharafutdinov \cite{CDS}
proved {\em local} filling minimality and boundary rigidity
for metrics with certain upper curvature bounds.

In this paper we settle the local version
of the minimal filling conjecture,
namely we prove the following theorem.

\begin{theorem}
\label{t-riemannian}
For every simple Riemannian metric $g_0$ on $D^n$
there is a neighborhood $\U$ of $g_0$ in the space
of all Riemannian metrics on $D^n$ (with the
$C^\infty$ topology) such that the following holds.
For all metrics $g,g'\in \U$ such that
$$
 d_{g'}(x,y) \ge d_g(x,y) \qquad\text{for all $x,y\in\pd D$} ,
$$
one has
$$
 \vol(D^n,g') \ge \vol(D^n,g) .
$$
\end{theorem}

The same result holds for Finsler metrics (see Theorem \ref{t-finsler} below),
and even in the Riemannian case the proof relies on Finsler geometry.
The plan of the proof is the following. First we show that the fact
that boundary distances determine the volume extends to Finsler
metrics. Then, for metrics $g$ and $g'$ as in Theorem~\ref{t-riemannian},
we construct a (non-reversible) Finsler metric $\phi$ on $D$ which induces
the same boundary distances as $g'$ and majorizes $g$ pointwise in~$D$.
It then follows that $\vol(D,g')=\vol(D,\phi)\ge\vol(D,g)$.
Finsler metrics help here because they are flexible:
unlike in the Riemannian case, it is easy to construct
perturbations of a Finsler metric that induce
a given perturbation of the boundary distance function.
The actual details of the proof are more complicated
than the above plan:
to work around non-smoothness of the boundary distance function
at the diagonal, we extend the metric to a slightly larger disc
and use distances in that larger disc rather than the original one,
see section \ref{sec-env}.

Let us proceed with definitions and formulations for the
Finslerian case. A \textit{Finsler manifold} is a smooth manifold
equipped with a Finsler metric.
A \textit{Finsler metric} on a smooth manifold $M$
is a continuous function $\phi\colon TM\to\R$ satisfying the following
conditions:

(1) $\phi(tv)=t\phi(v)$ for all $v\in TM$ and $t\ge 0$;

(2) $\phi$ is positive on $TM\setminus 0$;

(3) $\phi$ is smooth on $TM\setminus 0$;

(4) $\phi$ is strictly convex in the following sense:
for every $x\in M$, the function $\phi_x:=\phi^2|_{T_xM}$
has positive definite second derivatives on $T_xM\setminus\{0\}$.

Note that we do {\em not} require that 
$\phi(-v)=\phi(v)$, 
that is, non-reversible Finsler metrics
are allowed. Nevertheless we still referq
to functions $\phi_x$ as \textit{norms} on the fibers $T_xM$.
For a Finsler metric $\phi$, one naturally defines geodesics, lengths,
and a (non-symmetric) distance function
$d_\phi\colon M\times M\to\R_+$,
see e.g.\ \cite{BCS} 
for details.
We define the notion of a simple Finsler metric and
the boundary distance function $bd_\phi$
in the same way as in the Riemannian case.

For the notion of volume of a Finsler metric
we use the Holmes--Thompson definition \cite{HT},
reproduced here for the reader's convenience.
Let $M=(M,\phi)$ be a Finsler manifold.
Consider the co-tangent bundle $T^*M$ and
let $\phi^*\colon T^*M\to\R$ be the fiber-wise
dual norm to $\phi$. That is, for $x\in M$
and $\alpha\in T^*_xM$, one defines
$$
 \phi^*(\alpha) = \sup\{ \alpha(v) \mid v\in T_xM,\ \phi(v)=1 \} .
$$
Let $B^*M=B^*(M,\phi)$ 
the bundle of unit balls of $\phi^*$:
$$
 B^*M = \{ \alpha\in T^*M \mid \phi^*(\alpha) \le 1 \} .
$$
The \textit{Holmes--Thompson volume} of $M$,
that we denote by $\vol(M)$ or $\vol(M,\phi)$,
is defined by
$$
 \vol(M) = \frac1{\omega_n} \Vol_{can}(B^*M) 
$$
where $n=\dim M$, $\omega_n$ is the Euclidean volume
of the unit ball in $\R^n$, and $\Vol_{can}$ is the canonical
(symplectic) $2n$-dimensional volume on $T^*M$.
Clearly this definition yields the Riemannian volume
in the case when the Finsler metric is Riemannian.
Also notice that the volume is monotonous with respect
to the metric: if $\phi'\ge\phi$ pointwise, then
$B^*(M,\phi')\supset B^*(M,\phi)$
and therefore $\vol(M,\phi')\ge\vol(M,\phi)$.

\begin{theorem}
\label{t-finsler}
For every simple Finsler metric $\phi_0$ on $D^n$
there is a $C^\infty$ neighborhood $\U$ of $\phi_0$
such that the following holds.
For all Finsler metrics $\phi,\phi'\in \U$ such that
$ bd_{\phi'} \ge bd_\phi$ pointwise,
one has
$
 \vol(D^n,\phi') \ge \vol(D^n,\phi) 
$.
\end{theorem}

\begin{remark}
Similar Finsler volume comparison results
(in the reversible case)
were obtained recently by H.~Koehler \cite{K}. In particular,
Corollary 3.2(2) in \cite{K} asserts filling minimality
of $\phi$ among all simple metrics $\phi'$ such that
$bd_{\phi'}$ is sufficiently close to $bd_\phi$ in the
strong $C^2$ topology (in the complement of the diagonal
in $\pd D\times\pd D$). This result implies the assertion
of Theorem \ref{t-finsler} under an additional assumption
that $\phi$ and $\phi'$
(along with their derivatives up to a certain order)
agree on $\pd D$. 
\end{remark}

\begin{remark}
Finsler metrics are never boundary rigid as they admit non-isometric
perturbations preserving boundary distances.
(One possible construction is described
in section \ref{sec-env},
see the paragraph preceding Lemma \ref{l-vol-via-envelope}.)
Because of this, our Finslerian proof of Theorem~\ref{t-riemannian}
does not have immediate rigidity implications.
\end{remark}

Theorem \ref{t-riemannian} is a special case of
Theorem \ref{t-finsler}. The proof of Theorem \ref{t-finsler}
occupies sections \ref{sec-env} and \ref{sec-proof},
with Appendix~\ref{app} containing the proof of a
technical lemma.

In section \ref{sec-xray} we show how Theorem \ref{t-finsler}
implies the (well-known) injectivity of the geodesic ray transform
for a simple Finsler metric, see Corollary \ref{cor-xray}.
Loosely speaking, this injectivity means that a smooth
function on $D$ is uniquely determined by its integrals
over geodesics. More precisely, consider
a simple Finsler metric $\phi$ on $D=D^n$
and denote by $\Gamma_\phi$ the space of all maximal geodesics
of this metric (this space is a $(2n-2)$-dimensional smooth manifold
diffeomorphic to the complement of the diagonal in $\pd D\times\pd D$).
The \textit{geodesic ray transform} of $\phi$ is a map
$$
 I_\phi\colon C^\infty(D) \to C^\infty(\Gamma_\phi)
$$
defined by
$$
 I_\phi f(\gamma) = \int_{dom(\gamma)} f(\gamma(t))\,dt \qquad f\in C^\infty(D),\ \gamma\in\Gamma_\phi .
$$
(Here and everywhere in the paper the geodesics are
parametrized by arc length, that is, $\phi(\dot\gamma(t))\equiv 1$
for all $\gamma\in\Gamma_\phi$ and $t\in dom(\gamma)$.)

\begin{cor}\label{cor-xray}
If $n\ge 2$ and $\phi$ is a simple Finsler metric on $D^n$, then $I_\phi$ is injective.
\end{cor}

This result is not new; a more general theorem is
proved by Sharafutdinov \cite{sharaf} by analytic methods.
For Riemannian metrics, the injectivity of the geodesic
ray transform for a simple metric is proved
by Mukhometov \cite{Mu1,Mu2} and 
independently by Bernstein and Gerver \cite{BG1,BG2}.

Corollary \ref{cor-xray} easily follows from Theorem \ref{t-finsler} applied
to metrics conformal to $\phi$, see section \ref{sec-xray} for details.
The author believes that this new proof is more geometric and transparent
than the one in \cite{sharaf}.

\smallskip\textit{Acknowledgement}.
The author is grateful to V.~A.~Sharafutdinov and G.~Uhlmann
for useful discussion of the history of the geodesic ray transform problem.

\section{Enveloping functions}
\label{sec-env}

Let $\phi$ be a Finsler metric on $D=D^n$.
We denote by $UD$ and $U^*D$ the bundles
of unit spheres of $\phi$ and $\phi^*$,
respectively.
To emphasize the dependence on $\phi$ where needed,
we write $U(D,\phi)$ and $U^*(D,\phi)$.
We say that a smooth function $f\colon D\to\R$
is \textit{distance-like} (with respect to $\phi$)
if $\phi^*(d_xf)=1$ for all $x\in D$.

If $\phi$ is simple, then
the distance function $d_\phi(p,\cdot)$
of a point $p\in D$ satisfies this requirement
everywhere except at~$p$. To construct a
smooth distance-like function, consider a larger
disc $D^+\supset D$ and smoothly extend the metric $\phi$
to it. For a fixed extension, choosing $D^+$
sufficiently close to $D$ guarantees that $\phi$
is simple on $D^+$. Then for every $p\in \pd D^+$
the function $d_\phi(p,\cdot)|_D$ is smooth
and distance-like.

The (Finslerian) \textit{gradient} of a distance-like function
$f\colon D\to\R$ at $x\in D$,
denoted by $\grad_\phi f(x)$,
is the unique
tangent vector $v\in U_xD$ such that $d_xf(v)=1$.
In other words, $\grad_\phi f(x)$ is the Legendre
transform of the co-vector $d_xf$
with respect to the Lagrangian $\frac12\phi^2$.

For example, the gradient at $x$ of the distance function
$d_\phi(p,\cdot)$ of a simple Finsler metric $\phi$
is the velocity at the endpoint of the unique
minimizing geodesic from $p$ to~$x$.
We denote this velocity vector by $\overleftarrow{xp}$.

A \textit{gradient curve} of $f$ is a curve
$\gamma\colon[a,b]\to D$ such that
$\dot\gamma(t)=\grad_\phi(\gamma(t))$ 
for all $t\in[a,b]$. Clearly every gradient curve
of a distance-like function is a minimizing geodesic.

\begin{definition}
\label{d-env}
Fix a manifold $S$ diffeomorphic to $S^{n-1}$.
We say that a smooth function
$F\colon S\times D\to\R$ is an \textit{enveloping function}
for $\phi$ if the following two conditions are satisfied:

(i) for every $p\in S$, the function $F_p:=F(p,\cdot)$
is distance-like;

(ii) for every $x\in D$, the map $p\mapsto d_xF_p$
is a diffeomorphism from $S$ to  $U^*_xD$.
\end{definition}

We construct an enveloping function for a simple
metric $\phi$ as follows.
Consider a disc $D^+\supset D$ with metric $\phi$
extended as above. Identify $S$ with $\pd D^+$ and define
$$
 F(p,x) = F_p(x) = d_\phi(p,x)
$$
for all $p\in S=\pd D^+$ and $x\in D$. Then $F$ is an
enveloping function. Indeed, for a fixed $x\in D$,
the map $p\mapsto\grad_\phi F_p(x)=\overleftarrow{xp}$ is a diffeomorphism
between $S$ and $U_xD$ because the metric is simple.
Since the Legendre transform (from $U_xD$ to $U_x^*D$)
is also a diffeomorphism, so is the map 
$p\mapsto d_xF_p$ from $S$ to $U^*_xD$.

Some non-simple metrics admit enveloping functions as well.
For example, if $F\colon S\times D\to\R$ is an enveloping function
for a metric $\phi$ on $D$, and $D'\subset D$ is a sub-domain with
smooth (but not necessarily convex) boundary,
then $F|_{S\times D'}$ is an enveloping function
for $\phi|_{D'}$. However the existence of an enveloping
function implies that all geodesics are minimizing
(because they are gradient curves of distance-like functions)
and therefore have no conjugate points.

An enveloping function $F$ uniquely determines the metric.
Indeed, the unit sphere of the dual norm $\phi^*$
at every point $x\in D$ is the image of the map
$p\mapsto d_xF_p$ from $S$ to $T^*_xD$
and thus is determined by $F$. This unit sphere
determines the dual norm $\phi^*_x$ and the latter
determines the Finsler norm $\phi_x$.

Furthermore every $C^3$-small perturbation of $F$
yields an enveloping function of a Finsler metric.
Indeed, let $F$ be an enveloping function
for $\phi$ and $F'\colon S\times D\to\R$ be
$C^3$-close to $F$.
Define $F'_p=F'(p,\cdot)$ for every $p\in S$.
Then for every $x\in D$, the map
$p\mapsto d_xF'_p$ is $C^2$-close
to the similar map for $F$.
Therefore the image of this map is
a convex surface in $T^*M$,
and this surface is the unit sphere
of some norm $\phi^{\prime*}_x$.
The dual norm to $\phi^{\prime*}_x$
is a norm $\phi'_x$ on $T_xM$,
and the union of these norms over all $x\in D$
is a Finsler metric $\phi'$ for which
$F'$ is an enveloping function.

The boundary distance function of a simple metric $\phi$
is uniquely determined by the restriction $F|_{S\times\pd D}$
of an enveloping function $F$ to the boundary. Namely,
$$
 d_\phi(x,y) = \max_{p\in S} (F(p,y)-F(p,x))
$$
for $x,y\in\pd D$.
(The maximum is attained at a point $p\in\pd D$ 
such that the minimizing geodesic from $x$ to~$y$
is a gradient curve of $F_p$.)
Therefore every $C^3$-small perturbation
of $F$ in the interior of its domain
produces a perturbation of the metric preserving
the boundary distance function.

\begin{lemma}
\label{l-vol-via-envelope}
Let $F$ be an enveloping function for $\phi$.
Then $\vol(D,\phi)$ is uniquely determined
by the restriction $F|_{S\times\pd D}$.
\end{lemma}

Note that in this lemma we do not assume that $\phi$ is simple.
However, as explained above, existence of an enveloping function
implies simplicity of the metric provided that the boundary is strictly convex.
One could deduce the lemma from the fact that the boundary distance
function of a simple Finsler metric uniquely determines the volume,
but the author is not aware of a published proof of this fact.

\begin{proof}[Proof of Lemma \ref{l-vol-via-envelope}]
For every $x\in D$, define a differential $(n-1)$-form $V_x$
on $U^*_xD$ with values in $\Lambda^n T_x^*D$ as follows:
for $\alpha\in U^*_xD$ and $\eta\in\Lambda^{n-1} T_v(U^*_xD)$, let
$$
 V_x(\alpha)(\eta) = \frac1n \cdot \alpha\wedge i_*(\eta)
$$
where $i\colon U^*_xD\to T^*_xD$ is the inclusion map
and $i_*\colon \Lambda^{n-1}T_v(U^*_xD)\to \Lambda^{n-1}T_x^*D$
the induced mapping of $(n-1)$-vectors.
The integral of $V_x$ over $U^*_xD$ equals the volume of
the unit co-tangent ball $B^*_x$ (note that the volume
of a subset of $T^*_xD$ has an invariant meaning as
an element of $\Lambda^n T_x^*D$). Thus the definition
of the Holmes--Thompson volume can be written in the form
\begin{equation}
\label{e-vol-Vx}
 \vol(D,\phi) = \frac1{\omega_n} \int_D \left(x\mapsto\int_{U^*_xD} V_x\right)
\end{equation}
Let $F\colon S\times D\to\R$ be an enveloping function for $\phi$.
For $x\in D$, define a map $G_x\colon S\to U^*_xD$ by
$$
 G_x(p) = d_xF_p
$$
where $F_p=F(p,\cdot)$.
By the second requirement of Definition \ref{d-env},
$G_x$ is a diffeomorphism.
It induces a map
$G_x^*\colon \Lambda^{n-1} T^*U^*_x\to \Lambda^{n-1} T^*S$
which carries differential $(n-1)$-forms from $U^*_x$ to~$S$.
From \eqref{e-vol-Vx} we have
\begin{equation}
\label{e-vol-via-lambda}
 \vol(D,\phi) = \frac1{\omega_n} \int_D\left(x\mapsto\int_S G^*_xV_x\right)
 = \frac1{\omega_n} \int_S \lambda
 \end{equation}
where $\lambda$ is a differential $(n-1)$-form on $S$ given by
$$
 \lambda = \int_D(x\mapsto G^*_xV_x) .
$$
We are going to show that the value $\lambda(p)$ of $\lambda$
at every point $p\in S$
is determined by the restriction $F|_{S\times\pd D}$.
We have $\lambda(p) = \int_D \nu_p$ where
$\nu_p$ is a differential $(n-1)$-form on $M$
with values in $\Lambda^{n-1}T^*_pS$ defined by
$$
 \nu_p(x) = (G^*_x V_x)(p) .
$$
Then for tangent vectors $\xi_1,\dots,\xi_{n-1}\in T_pS$ we have
$$
\begin{aligned}
 \nu_p(x) (\xi_1\wedge\dots\wedge\xi_{n-1}) 
 &= V_x (dG_x(\xi_1)\wedge\dots\wedge dG_x(\xi_{n-1})) \\
 &= \frac1n \cdot G_x(p)\wedge dG_x(\xi_1)\wedge\dots\wedge dG_x(\xi_{n-1})
\end{aligned}
$$
where the vectors $dG_x(\xi_i)$ in the right-hand side
are regarded as elements of $T^*_xD$
(we implicitly use the inclusion of $T_v U^*D$ into $T^*_xD$
for $v=G_x(p)$ here).
Substituting the definition of $G_x$ yields
$$
 \nu_p(x) (\xi_1\wedge\dots\wedge\xi_{n-1})
 = \frac1n \cdot d_xF_p \wedge d_x F_{p,\xi_1}\wedge\dots\wedge d_x F_{p,\xi_{n-1}}
$$
where $F_{p,\xi_i}$ denotes the derivative of $F$ 
with respect to the first argument along
the tangent vector $\xi_i$ at $p$.
Thus
$$
 \lambda(p) (\xi_1\wedge\dots\wedge\xi_{n-1}) 
 = \frac1n\int_D dF_p \wedge dF_{p,\xi_1}\wedge\dots\wedge dF_{p,\xi_{n-1}} .
$$
The $n$-form under the integral in the right-hand side equals
the exterior derivative of the $(n-1)$-form
$F_p \cdot dF_{p,\xi_1}\wedge\dots\wedge dF_{p,\xi_{n-1}}$.
Therefore by Stokes' theorem
$$
 \lambda(p) (\xi_1\wedge\dots\wedge\xi_{n-1})
 = \frac1n\int_{\pd D} F_p \cdot dF_{p,\xi_1}\wedge\dots\wedge dF_{p,\xi_{n-1}} .
$$
The right-hand side is determined by the restriction $F|_{S\times\pd D}$,
hence so are $\lambda(p)$ and the volume $\vol(D,\phi)$ 
which is written in terms of $\lambda$ in \eqref{e-vol-via-lambda}.
\end{proof}

\begin{example}
\label{x-2d}
Consider the case of dimension $n=2$.
As the above proof shows, the volume can be expressed
in terms of an enveloping function $F=F(p,x)$ restricted
on $S\times\pd D$ as follows:
$$
 \vol(D,\phi) = \frac1{2\pi}\int_S\int_{\pd D} F(p,x)\cdot
 \frac{\pd^2 F}{\pd p\pd x}(p,x)\,dxdp
 = -\frac1{2\pi}\int_{S\times\pd D} \frac{\pd F}{\pd x} \cdot
 \frac{\pd F}{\pd p} \,dxdp .
$$
Let $f=bd_\phi$. We use the notation $x$ and $y$ for the
arguments of~$f$.
Since $f$ is a limit of enveloping functions (whose
first derivatives converge a.e.), we can use the same
formula with $f$ in place of $F$ (and $\pd D$ in place of $S$):
\begin{equation}
\label{e-vol2}
 \vol(D,\phi) = -\frac1{2\pi}\int_{\pd D\times\pd D} \frac{\pd f}{\pd x} \cdot
 \frac{\pd f}{\pd y} \,dxdy .
\end{equation}
This is an explicit formula for the volume in terms of the boundary
distance function in dimension~2.

Identify $\pd D$ with the standard circle (of length $2\pi$).
Let us further restrict ourselves to the case when $f$
is symmetric and invariant under rotations
of the circle. Then $\rho$ is determined by a function
$f_0=f(x_0,\cdot)$ in one variable which ranges over
a half-circle. Identifying the half-circle with the segment
$[0,\pi]$, we can rewrite \eqref{e-vol2} as follows:
$$
 \vol(D,\phi) = 2 \int_0^\pi (\dot f_0)^2
$$
It is easy to see that this formula is not monotone
with respect to $f_0$, even within the class of
increasing functions satisfying the triangle inequality.
However it is monotone within the class of concave functions,
and any rotation-invariant boundary distance function
is concave.
\end{example}

\section{Proof of Theorem \ref{t-finsler}}
\label{sec-proof}

Fix a simple Finsler metric $\phi_0$ on $D$
and let $\phi,\phi'$ be Finsler metrics
$C^\infty$-close to $\phi_0$
and such that $bd_{\phi'}\ge bd_\phi$.
Since simplicity of a metric is an open condition,
we may assume that $\phi$ and $\phi'$ are simple.
As explained in the previous section, extend $\phi$ and $\phi'$
to a larger disc $D^+\supset D$ and define enveloping
functions
$$
F_0,F,F'\colon S\times D\to\R
$$
for $\phi_0$, $\phi$ and $\phi'$
respectively, by
$$
 F_0(p,x) = d_{\phi_0}(p,x), \qquad F(p,x) = d_\phi(p,x), \qquad F'(p,x) = d_{\phi'}(p,x)
$$
for $p\in S=\pd D^+$ and $x\in D$.
The extension of the metrics can be chosen so that
$F$ and $F'$ are $C^\infty$-close to $F_0$.

We are going to construct a function $F''\colon S\times D\to\R$
which is also $C^\infty$-close to $F_0$ (so that there is a metric
$\phi''$ for which $F''$ is an enveloping function)
and satisfies the following conditions:

(i) $F''|_{S\times\pd D}=F'|_{S\times\pd D}$;

(ii) the resulting metric $\phi''$ satisfies $\phi''\ge\phi$ pointwise.

The first condition implies that $\vol(D,\phi'')=\vol(D,\phi')$
by Lemma \ref{l-vol-via-envelope}.
The second one implies that $\vol(D,\phi'')\ge\vol(D,\phi)$.
Thus the existence of such $F''$ implies the desired
inequality $\vol(D,\phi')\ge\vol(D,\phi)$.

First we construct an auxiliary
map $\Psi:U(D,\phi')\to U(D,\phi)$ between the
unit tangent bundles of our metrics.
This map sends every trajectory of the geodesic flow
of $\phi'$ to the corresponding trajectory of
the geodesic flow of $\phi$ (with a linear change of time),
where ``corresponding''
means that the associated geodesics of the two metrics
connect the same pair of boundary points of~$D$.
More precisely,
let $x\in D$ and $v\in U_x(D,\phi')$,
and let $\gamma'\colon[0,T']\to D$ be the 
maximal $\phi'$-geodesic through $v$,
that is, $\gamma'(t')=x$ and $\dot\gamma'(t')=v$ for some $t'\in[0,T']$.
Let $\gamma\colon[0,T]\to D$ be the $\phi$-geodesic 
with the same endpoints on $\pd D$: $\gamma(0)=\gamma'(0)$ and $\gamma(T)=\gamma'(T')$.
If $T'\ne 0$, define
$
 \Psi(v) = \dot\gamma(t'T/T')
$,
and if $T'=0$ (or, equivalently, if $x\in\pd D$
and $v$ is tangent to $\pd D$), let $\Psi(v)$ be
the $\phi$-unit vector positively proportional to~$v$.

Clearly the map $\Psi$ defined this way is a homeomorphism
between $U(D,\phi')$ and $U(D,\phi)$.
Moreover $\Psi$ is a diffeomorphism and $\Psi$ goes
to the identity as $\phi,\phi'\to\phi_0$ in $C^\infty$,
see Proposition~\ref{p-psi} and Remark \ref{r-psi} in Appendix \ref{app}.

For every $p\in S$, define a map $G'_p\colon D\to U(D,\phi')$
by $G'_p(x)=\grad_{\phi'} F'_p(x)$.
Note that $\pi\circ G'_p=id_D$
where $\pi\colon T^*D\to D$ denotes the bundle projection.
Since $\Psi$ is close to the identity,
the map
$$
H_p:=\pi\circ\Psi\circ G'_p\colon D\to D
$$
is also $C^\infty$-close to the identity
and hence is a diffeomorphism.
Define a function $F''_p:D\to\R$ by
$F''_p=F'_p\circ H_p^{-1}$
and a function  $F''\colon S\times D\to\R$ by
$F''(p,x)=F''_p(x)$ for all $p\in S$, $x\in D$.
Note that $F''$ depends continuously on $\phi$ and $\phi'$
and $F''=F_0$ if $\phi=\phi'=\phi_0$, 
hence $F''\to F_0$ as $\phi,\phi'\to\phi_0$.
Therefore, if $\phi$ and $\phi'$ is sufficiently close to $\phi_0$,
there is a metric $\phi''$ such that $F''$
is an enveloping function for it.

It remains to verify that $F''$ and $\phi''$ satisfy
the above conditions (i) and (ii).
By construction, we have $H_p|_{\pd D}=id_{\pd D}$
for every $p\in S$, hence $F''_p|_{\pd D}=F'_p|_{\pd D}$
and therefore $F''$ satisfies (i).
To verify (ii), fix $x\in D\setminus\pd D$ and recall that
the unit sphere of the dual norm $\phi_x''^*$
is parametrized by the family $\{d_x F''_p\}_{p\in S}$.
Therefore
\begin{equation}
\label{e-phi-dF}
 \phi''(v) = \sup_{p\in S} \{ d_xF''_p(v) \}
\end{equation}
for every $v\in T_xD$.
Let $v\in U_x(D,\phi)$ and $\gamma:[0,T]\to D$
be the maximal $\phi$-geodesic through $v$,
that is,
$x=\gamma(t_0)$ and $v=\dot\gamma(t_0)$
for some $t_0\in[0,T]$.
Let $\gamma'\colon[0,T']\to D$ be the $\phi'$-geodesic
connecting the same endpoints on $\pd D$.
Then, by the definition of $\Psi$, we have
$$
 \Psi(\dot\gamma'(t')) = \dot\gamma(t' T/T')
$$
for all $t'\in[0,T']$.
Since $F'$ is an enveloping function for $\phi'$,
the $\phi'$-geodesic $\gamma'$ is a gradient curve
of the function $F'_p:=F'(p,\cdot)$ for some $p\in S$,
that is,
$$
 \dot\gamma'(t') = G'_p(\gamma'(t'))
$$
for all $t'\in[0,T']$. Therefore
$$
 H_p(\gamma'(t')) = \pi\circ\Psi\circ G'_p(\gamma'(t'))
 = \pi\circ\Psi(\dot\gamma'(t'))
 = \gamma(t'T/T')
$$
for all $t'\in[0,T']$ or, equivalently,
$$
 H_p^{-1}(\gamma(t)) = \gamma'(tT'/T)
$$
for all $t\in[0,T]$.
Thus
$$
 F''_p(\gamma(t)) = F'_p\circ H_p^{-1} = F'_p(\gamma'(tT'/T)) ,
$$
for all $t\in[0,T]$, hence
$$
 d_x F''_p (v) = \frac d{dt} F''_p(\gamma(t))\Big|_{t=t_0} 
 = \frac{T'}T\cdot \frac d{dt'} F'_p(\gamma'(t')) \Big|_{t'=t_0T'/T}
 = \frac{T'}T
$$
since $\gamma'$ is a $\phi'$-gradient curve of $F'$
and therefore $\frac d{dt'} F'_p(\gamma'(t'))\equiv 1$.
Recall that $T=d_\phi(a,b)$ and $T'=d_{\phi'}(a,b)$
where $a=\gamma(0)=\gamma'(0)$
and $b=\gamma(T)=\gamma'(T')$,
and $d_{\phi'}(a,b)\ge d_\phi(a,b)$ since $bd_{\phi'}\ge bd_\phi$.
Therefore $d_x F''_p (v)=T'/T\ge 1$.

This and \eqref{e-phi-dF} imply that $\phi''(v)\ge 1$.
Since $v$ is an arbitrary $\phi$-unit tangent vector
at $x$, and $x$ is an arbitrary interior point of $D$,
it follows that $\phi''\ge\phi$.
This completes the proof of (ii) and hence of Theorem \ref{t-finsler}.

\section{Injectivity of geodesic ray transform}
\label{sec-xray}

The goal of this section is to deduce Corollary \ref{cor-xray}
from Theorem \ref{t-finsler}.
Let $\phi$ be a simple Finsler metric on $D=D^n$
and $I_\phi\colon C^\infty(D)\to C^\infty(\Gamma_\phi)$
its geodesic ray transform.
Let $f\in C^\infty(D)$ be such that $I_\phi f=0$;
we are to show that $f=0$.

For a small $\ep>0$, define a Finsler metric $\phi_\ep$ on $D$ by
$$
\phi_\ep(v)=(1+\ep f(\pi(v)))\cdot\phi(v), \qquad v\in TD
$$
where $\pi:TD\to D$ is the bundle projection.
Let $p,q\in\pd D$ and $\gamma:[0,T]\to D$ the
$\phi$-geodesic connecting $p$ to $q$.
Its $\phi_\ep$-length $L_{\phi_\ep}(\gamma)$ satisfies
$$
 L_{\phi_\ep}(\gamma) = \int_0^T (1+\ep f(\gamma(t))\, dt
 = T + \ep I_\phi f(\gamma) = T
$$
since $I_\phi f=0$. Since $\phi$ is simple, we have $T=L_\phi(\gamma)=d_\phi(p,q)$.
Therefore
$$
 d_{\phi_\ep}(p,q) \le L_{\phi_\ep}(\gamma) =d_\phi(p,q).
$$
Thus $bd_{\phi_\ep}\le bd_\phi$.
Since $\phi_\ep\to\phi$ in $C^\infty$ as $\ep\to 0$, Theorem \ref{t-finsler}
applies and we conclude that $\vol(D,\phi_\ep)\le\vol(D,\phi)$
for a sufficiently small $\ep$.
On the other hand,
$$
 \vol(D,\phi_\ep) = \int_D (1+\ep f)^n\, d\vol_\phi
$$
where $\vol_\phi$ is the volume form of $\phi$, therefore
$$
 \int_D (1+\ep f)^n\, d\vol_\phi \le \vol(D,\phi) .
$$
The same argument applied to $-f$ in place of $f$ yields that
$$
 \int_D (1-\ep f)^n\, d\vol_\phi \le \vol(D,\phi) .
$$
Summing these two inequalities we obtain
$$
 \int_D \big((1+\ep f)^n+(1-\ep f)^n - 2)\, d\vol_\phi \le 0 .
$$
This and the trivial inequality $(1+\ep f)^n+(1-\ep f)^n-2 \ge \ep^2 f^2$
for $n\ge 2$ imply that $\int_D f^2\,d\vol_\phi\le 0$, hence $f=0$.
This completes the proof of Corollary \ref{cor-xray}.

\appendix
\section{Smoothness of $\Psi$}
\label{app}

The goal of this appendix is to prove the following technical fact
used in the proof of Theorem \ref{t-finsler}.

\begin{prop}
\label{p-psi}
Let $\phi$ and $\phi'$ be arbitrary simple Finsler metrics on $D=D^n$
and let a map $\Psi\colon U(D,\phi')\to U(D,\phi)$
be defined as in section \ref{sec-proof}.
Then $\Psi$ is a $C^\infty$ diffeomorphism
and it depends smoothly on $\phi$ and $\phi'$.
\end{prop}

\begin{remark}
Technically, the domain $U(D,\phi')$ of $\Psi$ is a variable
(i.e., depending on $\phi'$) submanifold of $TD$.
To formalize the notion of smooth dependence on $\phi'$
in Proposition \ref{p-psi},
one can identify all unit tangent bundles by means of the
fiber-wise radial projection.
\end{remark}

\begin{remark}
\label{r-psi}
Obviously $\Psi$ is the identity in the case when $\phi'=\phi$.
Therefore the smooth dependence on the metrics in
Proposition \ref{p-psi} implies
that $\Psi$ goes to the identity (in $C^\infty$)
as $\phi,\phi'\to\phi_0$.
\end{remark}

\begin{proof}[Proof of Proposition \ref{p-psi}]
We write $\Psi=\Psi_{\phi',\phi}$ to emphasize the dependence on the metrics.
We identify $D$ with the standard unit ball in $\R^n$ and
denote by $\phi_e$ the standard Euclidean metric
(regarded as a Finsler metric on $D$).
We will show that, for any simple metric $\phi$,
the map $\Psi_{\phi,\phi_e}$ is a diffeomorphism
and it depends smoothly on $\phi$.
Proposition \ref{p-psi} follows from this special case and the
trivial identity
$\Psi_{\phi',\phi} = \Psi_{\phi_e,\phi}\circ \Psi_{\phi',\phi_e}$.

For $u\in UD:=U(D,\phi)$, let $\gamma_u$ denote the maximal $\phi$-geodesic defined
by the initial data $\dot\gamma(0)=u$ and let $[\tau^-(u),\tau^+(u)]$
be the domain of $\gamma_u$.
Let $\ell(u)=\tau^+(u)-\tau^-(u)$ be the length of the geodesic $\gamma_u$,
$p^\pm(u)=\gamma_u(\tau^\pm(u))$ its endpoints on the boundary,
and $\tau(u)=\frac12(\tau^+(u)+\tau^-(u))$ the parameter of its midpoint.
Define
$$
\begin{aligned}
 \lambda(u) &= \frac{|p^+(u)-p^-(u)|}{\ell(u)} \in \R, \\
 p(u) &= \frac{|p^+(u)+p^-(u)|}2 \in D\subset\R^n, \\
 w(u) &= \frac{p^+(u)-p^-(u)}{|p^+(u)-p^-(u)|} \in S^{n-1} \\
\end{aligned}
$$
if $u$ is not tangent to $\pd D$.
If $u$ is tangent to $\pd D$, 
we have $\ell(u)=0$ and $p(u)=p^+(u)=p^-(u)$
and extend $\lambda$ and $w$ by continuity:
$\lambda(u)=|u|$, $w(u)=u/|u|$.
Then the map $\Psi=\Psi_{\phi,\phi_e}$ can be written
in the form
$$
 \Psi(u) = \big( p(u)-\tau(u)\lambda(u)w(u), w(u) \big) \in D\times S^{n-1} = U(D,\phi_e) .
$$
We are going to show that the functions $\tau$, $\lambda$, $p$ and $w$
(and hence $\Psi$) are smooth on $UD$.
Let $V\subset UD$ denote the set of $\phi$-unit vectors tangent to $\pd D$. 
This is a $(2n-3)$-dimensional submanifold of the boundary $\pd UD$.
The above functions are obviously smooth away from $V$,
so we only need to prove their smoothness at $V$.

Extend the metric $\phi$ to an open ball $D^+\supset D$.
Since $\pd D$ is strictly convex with respect to $\phi$, the trajectories
of the geodesic flow are nowhere tangent to~$V$.
Therefore in a neighborhood of $V$ in $UD^+$
there exists a coordinate system $(t,y,v)$
where $t,y\in\R$ and $v\in V$, such that
the $t$-lines are trajectories of the geodesic flow
and the points of $V$ have coordinates $(t,y,v)$ with $t=y=0$.

Observe that trajectories of the geodesic flow are tangent to $\pd UD$
at~$V$. Hence, by the implicit function
theorem, in a suitable neighborhood of $V$ the set $\pd UD$ is represented
by a coordinate equation $y=h_v(t)=h(t,v)$ where $h\in C^\infty(\R\times V)$,
and $h_v(0)=h_v'(0)=0$. Since $\pd D$ is strictly convex, we have $h_v''(0)\ne 0$
and we may assume that $h_v''(0)>0$ (changing coordinate $y$ to $-y$ if
necessary). Then the set $UD\subset UD^+$ is locally the set of solutions
of the inequality $y\ge h_v(t)$ in our coordinates.
We need the following standard lemma.

\begin{lemma}
\label{l-anal}
Let $f\in C^\infty(\R)$. Then

1. If $f(0)=0$, then there exists $g\in C^\infty(\R)$
such that $f(x)=xg(x)$ for all $x\in\R$.

2. If $f$ is even, i.e.\ $f(x)=f(-x)$ for all $x\in\R$,
then there exists $g\in C^\infty(\R)$
such that $f(x)=g(x^2)$ for all $x\in\R$.

In both cases, $g$ depends smoothly on $f$.
\qed
\end{lemma}

Since $h_v(0)=h_v'(0)=0$, we can apply the first part of
Lemma \ref{l-anal} to $h_v$ twice and
conclude that $h_v(t)=t^2g_v(t)$ for some smooth
function $g_v$.
Observe that $g_v(0) = \frac12 h_v''(0)>0$,
so $g_v$ is positive in a neighborhood of~0.
In this neighborhood we have
$h_v(t)=f_v(t)^2$ where $f_v(t)=t\sqrt{g_v(t)}$,
so $f_v\in C^\infty$ and $f'(0)>0$.
Let $u\in UD$ have coordinates $(t,y,v)$
where $t$ and $y$ are close to~0.
Then the values $\tau^\pm(u)$ 
satisfy the equation
$$
 h_v(t+\tau^\pm(u)) = y
$$
since the point with coordinates $(t+\tau^\pm(u),y,v)$
belongs to $\pd UD$. Therefore
$$
 f_v(t+\tau^\pm(u)) = \pm\sqrt y .
$$
Since $f_v'(0)>0$, $f_v$ is invertible near 0 and
the above equation implies that
\begin{equation}
\label{e-tau}
\begin{aligned}
 \tau^+(u) &= f_v^{-1}(\sqrt y) - t , \\
 \tau^-(u) &= f_v^{-1}(-\sqrt y) - t ,
\end{aligned}
\end{equation}
hence
$$
 \tau(u) = \tfrac12 \big(f_v^{-1}(\sqrt y)+f_v^{-1}(-\sqrt y)\big) - t .
$$
The function $x\mapsto f_v^{-1}(x)+f_v^{-1}(-x)$
is defined in a neighborhood of~0, smooth and even.
Hence, by the second part of Lemma \ref{l-anal},
this function has a form $x\mapsto\theta_v(x^2)$ where $\theta_v$ is a smooth function
(which depends smoothly on $v$). Then
$$
 \tau(u) = \tfrac12 \theta_v(y) - t
$$
and therefore $\tau$ is smooth in a neighborhood of~$V$.

To prove the smoothness of the map $p\colon UD\to D$, observe that
\begin{equation}
\label{e-p}
\begin{aligned}
  p^+(u) &= \pi(t+\tau^+(u),y,v) = \pi(f_v^{-1}(\sqrt y),y,v), \\
  p^-(u) &= \pi(t+\tau^-(u),y,v) = \pi(f_v^{-1}(-\sqrt y),y,v) \\
\end{aligned}
\end{equation}
where $\pi$ is the bundle projection $UD\to D$ represented as a function of coordinates.
Hence
$$
 p(u) = \tfrac12 \big( \pi(f_v^{-1}(\sqrt y),y,v) + \pi(f_v^{-1}(-\sqrt y),y,v)\big) .
$$
Since the map 
$$
(x,y,v)\mapsto \pi(f_v^{-1}(x),y,v) + \pi(f_v^{-1}(-x),y,v)
$$
is smooth and even with respect to the $x$ variable,
by the second part of Lemma \ref{l-anal} it has the form
$(x,y,v) \mapsto \alpha(x^2,y,v)$ where $\alpha$ is a smooth function.
Then $p(u)=\frac12\alpha(y,y,v)$,
hence $p$ is smooth.

To handle the functions $\lambda$ and $w$, we need the following lemma.

\begin{lemma}
\label{l-anal2}
For every $f\in C^\infty(\R)$ there exists $g\in C^\infty(\R)$
such that
$$
 f(\sqrt y)-f(-\sqrt y) = \sqrt y\cdot g(y)
$$
for all $y\ge 0$.
Furthermore, $g$ depends smoothly on $f$.
\end{lemma}

\begin{proof}
Define $f_1(x) = f(x)-f(-x)$ for all $x\in\R$.
Since $f_1$ is a smooth odd function, by the first part
of Lemma \ref{l-anal} it can be written in the form
$f_1(x)=xf_2(x)$ where $f_2$ is an even smooth function
(depending smoothly on $f$).
By the second part of Lemma \ref{l-anal},
$f_2$ can be written in the form $f_2(x)=g(x^2)$
where $g$ is a smooth function
(depending smoothly on $f_2$ and hence on $f$).
Thus $f(x)=xg(x^2)$.
Substituting $x=\sqrt y$ yields the result.
\end{proof}

Lemma \ref{l-anal2}, \eqref{e-tau} and \eqref{e-p} imply that
$$
 \ell(u) = \tau^+(u) - \tau^-(u) = f_v^{-1}(\sqrt y) - f_v^{-1}(-\sqrt y)
 = \sqrt y \tilde\ell(y,v)
$$
and
$$
 p^+(u)-p^-(u) = \pi(f_v^{-1}(\sqrt y),y,v)-\pi(f_v^{-1}(-\sqrt y),y,v)
 = \sqrt y \tilde w(y,v)
$$
where $\tilde\ell$ and $\tilde w$ are smooth functions.
Since $f_v^{-1}$ has nonzero derivative at 0,
we have $\tilde\ell(0,v)\ne 0$ and $\tilde w(0,v)\ne 0$.
Now we have
$$
 \lambda(u) = \frac{|p^+(u)-p^-(u)|}{\ell(u)} = \frac{|\tilde w(y,v)|}{\tilde\ell(y,v)}
$$
and
$$
 w(u) = \frac{p^+(u)-p^-(u)}{|p^+(u)-p^-(u)|} = \frac{\tilde w(y,v)}{|\tilde w(y,v)|}
$$
where the right-hand sides are smooth in a neighborhood of the set $\{y=0\}$.
Thus $\lambda$ and $w$ are smooth in a neighborhood of $V$
and therefore the map $\Psi=\Psi_{\phi,\phi_e}$ is smooth on $UD$.

It remains to verify that $\Psi$ is a diffeomorphism depending smoothly
on the metric~$\phi$. It is easy to see from the definition that
$\Psi^{-1}=\Psi_{\phi_e,\phi}$ is smooth away from the set $\Psi(V)$
of unit vectors tangent to $\pd D$. The strict convexity of the boundary
easily implies that the derivative of $\Psi$ is non-degenerate at any point of~$V$.
Therefore $\Psi$ is a diffeomorphism. To verify the smooth dependence on $\phi$,
observe that the extension of $\phi$ to $D^+$ and the coordinates $(t,y,v)$
can be constructed in such a way that they depend smoothly on $\phi$.
Then all the smooth functions constructed throughout the proof depend
smoothly on $\phi$ and hence so does the map $\Psi$.
\end{proof}

\bibliographystyle{plain}

\end{document}